\documentclass{article}
\usepackage{amsfonts}
\usepackage[dvips]{graphicx}
\usepackage{color} 
\begin{document}

\title{Variational Method \\for Optimal Multimaterial Composites and Optimal Design}
\maketitle
\setcounter{page}{1}
\abstract{The paper outlines novel variational technique for finding microstructures of optimal multimaterial composites, bounds of composites properties, and multimaterial optimal designs. The translation method that is used for the exact two-material bounds is complemented by additional pointwise inequalities on  stresses in materials within an optimal composites.  The method leads to exact multimaterial bounds and provides a hint for optimal structures that may be multi-rank laminates or, for isotropic composites,  ``wheel assemblages". The Lagrangian of the formulated nonconvex multiwell variational problem is equal to the energy of the best adapted to the loading microstructures plus the cost of the used materials; the technique improves both the lower and upper bounds for the quasiconvex envelope of that Lagrangian. 
The problem of 2d elastic composites is described in some details; on particular, the isotropic component of the quasiconvex envelope of three-well Lagrangian for elastic energy is computed. The obtained results are applied  for computing of optimal multimaterial elastic designs;  an example of such a design is demonstrated. Finally, the optimal ``wheel assemblages"  are generalized and novel types of exotic microstructures with unusual properties are described.}

\section{Introduction}

Modern technologies of microfabrication and  3d printing allow for a huge variety of structures to be
manufactured for roughly the same price. Naturally, the material scientists want to know  what is ``the best'' structure, or how composites microstructures can be optimized; these questions are also related to metamaterials that utilize various extreme properties.   A close problem is the range of improvement of overall composite properties  that can be achieved by varying the structure.  There is no boundary between optimal design and an optimal composite material, which is also a structure at the microlevel: optimal designs are made from optimal composites. 
  So far, the vast majority of related results deals with two-material composites because of theoretical limitations. Meanwhile,  numerous applications call for optimal design of multimaterial composites, or even of porous composites made of two materials and void. Such designs are crucial  for  multi-physics applications, i.e. piezo-magnetic  and electromagnetic devices, in metamaterials and adaptive structures. 

\begin{figure}
 \centering
\begin{center}
 \includegraphics[width=2.5in,height=1.8in]{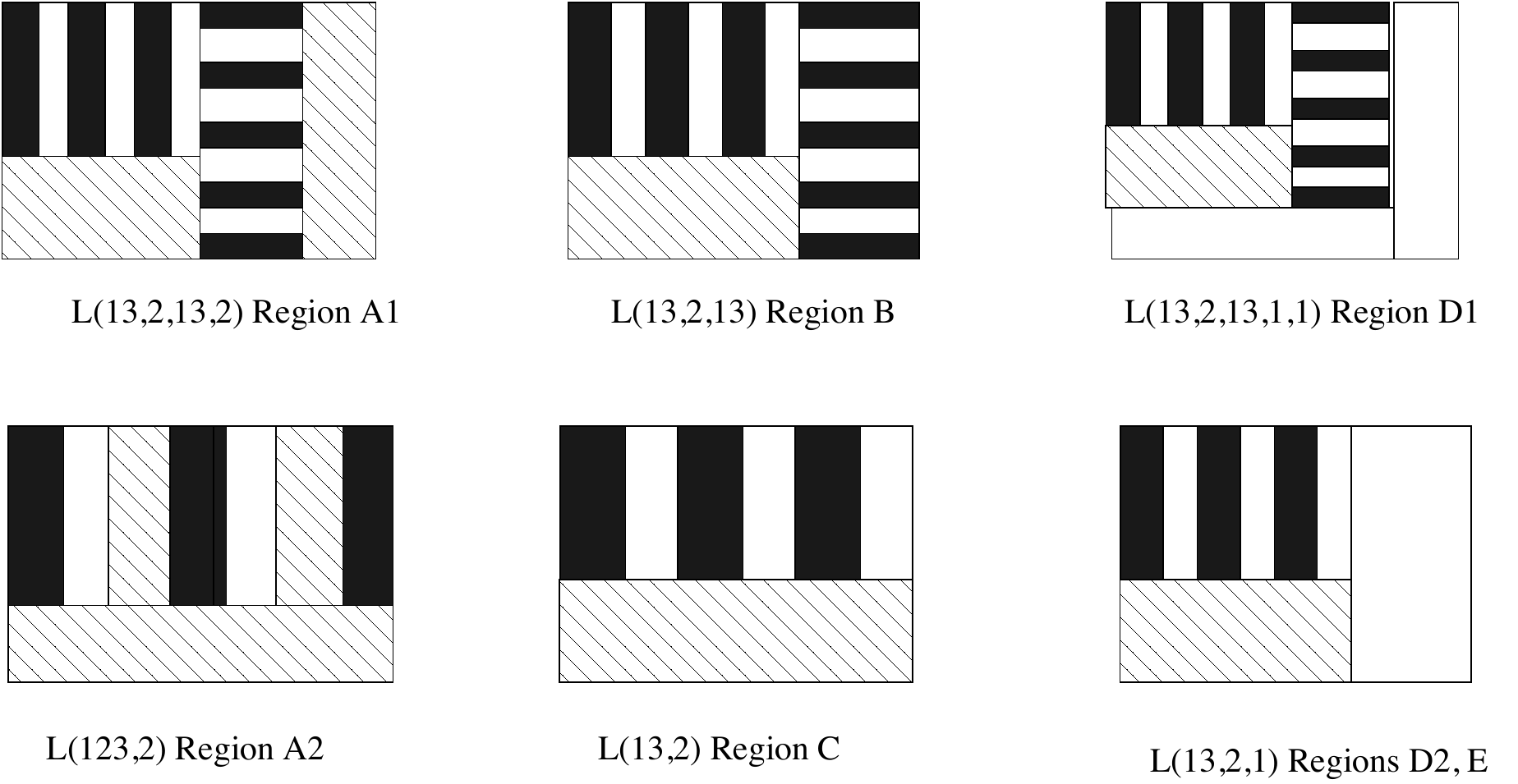} ~~  \includegraphics[width=3.7in,height=2.2in]{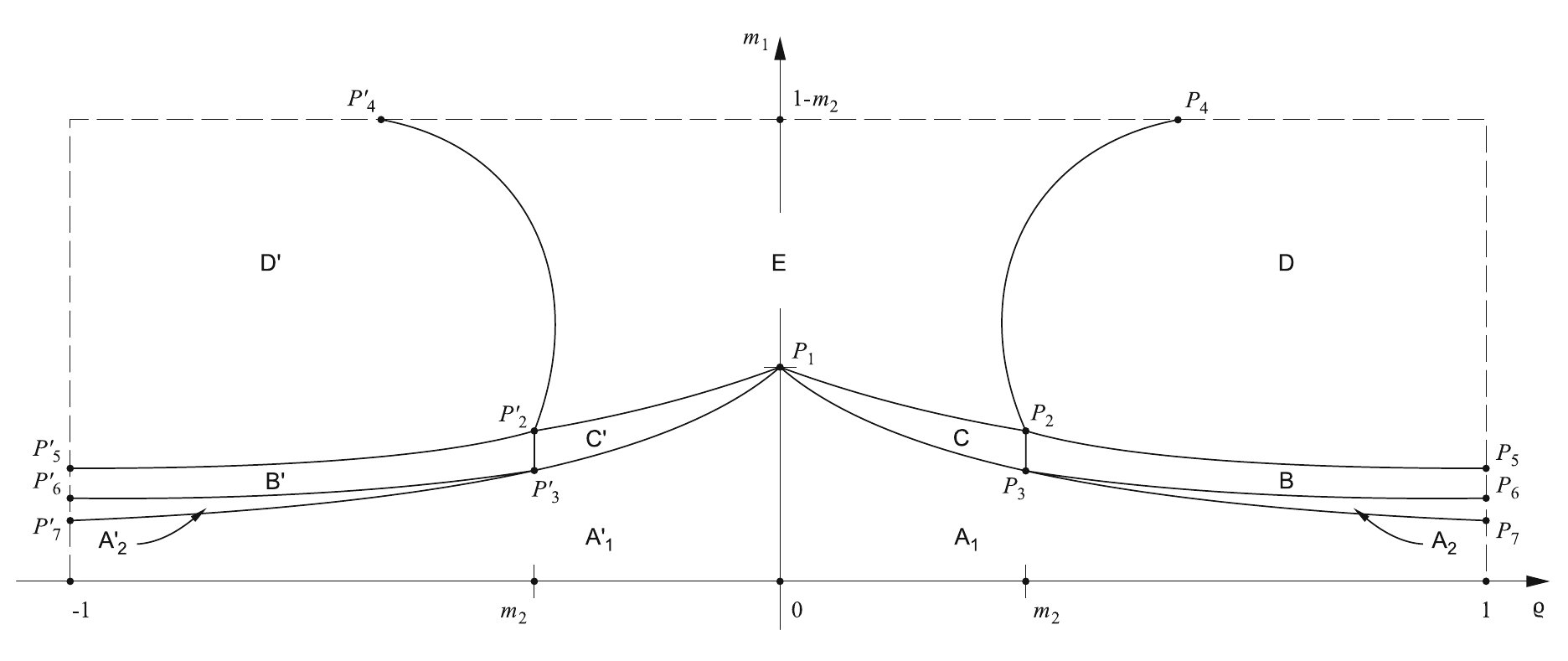} 
 \\
\caption{
{\it Left}: Cartoon of optimal multi-rank laminates that minimize elastic energy (compliance) of a three-material composite, see \cite{cherkyuan,cherkgz}. The parameters and the types A-E of the structures depend on the volume fractions and the ratio $p$  of  eigenvalues of the applied external stress $\sigma_0$. Black fields denote void (an infinite compliance,  $\kappa_3=\infty$), striped fields denote a material of intermediate compliance $\kappa_2$ and white fields demote the stiffest material $\kappa_1$, $\kappa_1 < \kappa_2$. The notation L(13,2,13)  shows the order of laminating as follows: materials $\kappa_1$ and $\kappa_3$ are laminated first, than they are laminated with material $\kappa_2$ in a orthogonal direction, then again laminated in an orthogonal direction with $\kappa_1$-$ \kappa_3$ laminate.   
\newline {\it Right}: 
Regions of optimality of the structures A-E  in dependence of the volume fraction $m_1$ of the best material (vertical axis) and $p$  (horizontal axis) \cite{cherkgz}. Volume fraction of $\kappa_2$ is fixed. The right vertical line corresponds to uniform pressure, the center vertical line corresponds to uniaxial load, the left vertical line corresponds to pure shear load, see below, Section \ref{newbound}.}
\label{opt-str}
\end{center}
\end{figure}

\begin{figure}
 \centering
\begin{center}
  \includegraphics[width=1.6in,height=.9in]{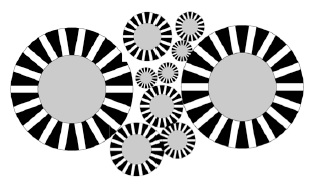}  ~~~~
  \includegraphics[width=.8in,height=.8in]{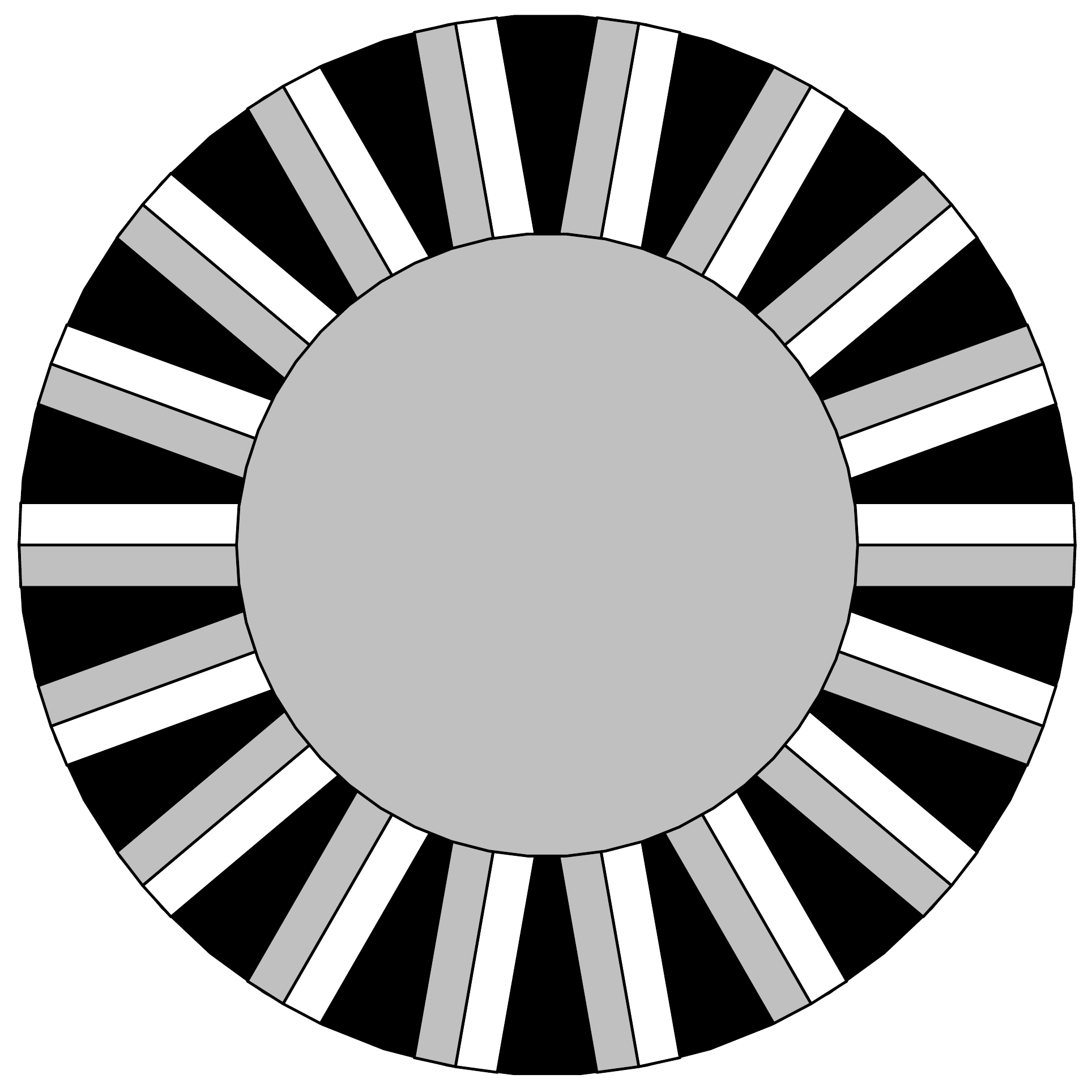} ~~ \includegraphics[width=.8in,height=.8in]{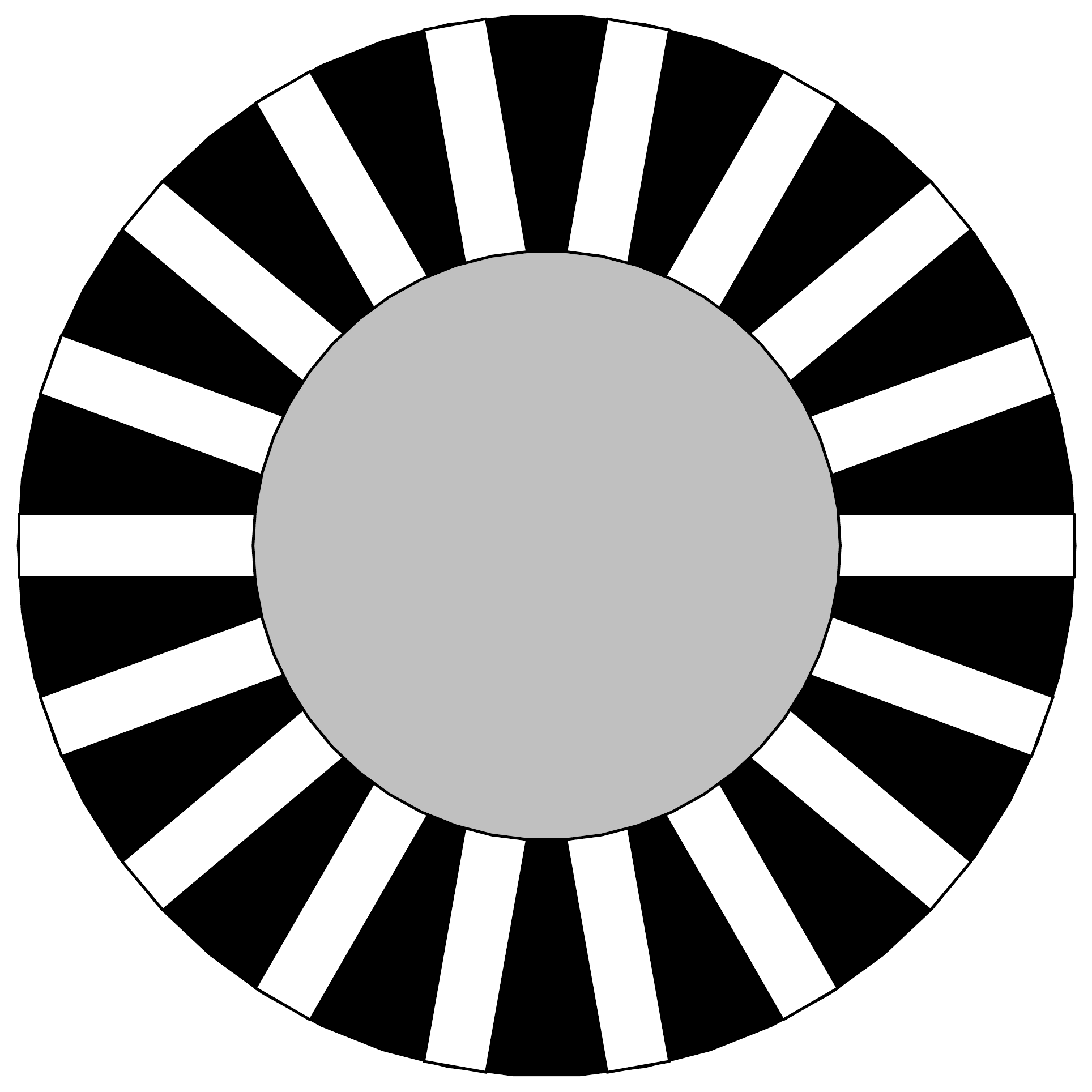}
 ~~ \includegraphics[width=.8in,height=.8in]{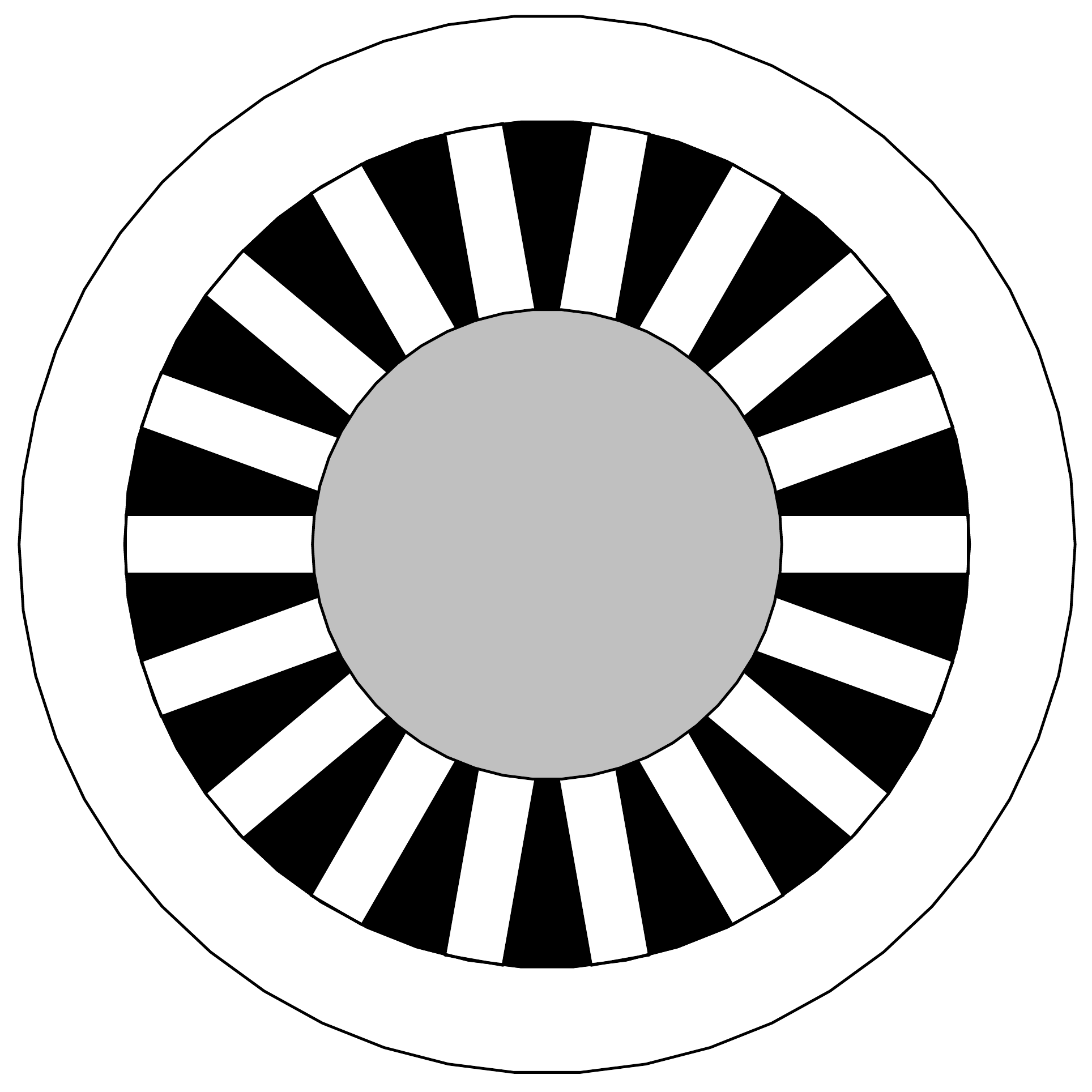} 
\caption{ An alternative optimal wheel-type assemblage for optimal isotropic microstructures (left) and elements W1, W2, W3 of the assemblage in dependence on volume fractions, see \cite{wheel}. The black field here denotes void, the grey denotes material $\kappa_2$, and the white denotes  material $\kappa_1$. The increase of the fraction $m_1$ of $\kappa_1$ (from left to right) leads to two topological transitions. The bulk modulus of the assemblages is equal to the bulk modulus of corresponding laminates in Figure \ref{opt-str}  made from the same materials taken in the same proportions. }
\label{wheels}
\end{center}
\end{figure}

Optimal microstructures of  multimaterial composites differ drastically from two-material ones. The latter have a steady and intuitively expected topology: 
 a strong material always surrounds weak inclusions, as in Hashin-Shtrikman coated circles and second-rank laminates which may degenerate to simple laminates.  
In contrast, optimal three-material structures  \cite{cherkyuan,cherkgz,wheel} (Figures \ref{opt-str} and \ref{wheels}) show a large variety of patterns and the optimal topology depends on the volume fractions. Optimal structures are diverse; they may or may not contain a strong envelope, and they may contain ``hubs" of intermediate material connected by anisotropic ``pathways" - laminates from the strong and weak materials, envelopes, and other configurations that reveal a geometrical essence of optimality. These structures are not unique, as shown in Figures \ref{opt-str} and \ref{wheels}.  

Obviously,  the methods  for finding them differ from the already developed methods used for optimal two-material structures. 
In this paper, we outline methods  for determination of multimaterial optimal elastic composites and designs from them. The exposition is partially based on results obtained in \cite{cherk09,cherkyuan,wheel,cherkpruss,cherkgz,bcd}

\section{Problems about optimal composites}
 \paragraph{The problem} The problem of the structure of optimal multimaterial composite has been studied for several decades. The bounds for multimaterial composites problem have been investigated  starting from the papers by Hashin and Shrikman \cite{HS1}, Milton \cite{milton82}, Lurie \& Cherkaev \cite{lc-multi},  Kohn \& Milton \cite{kohnmilton}. In 1995, Nesi \cite{nesi} suggested bounds for multimaterial mixtures that are better than Hashin-Shtrikman bounds; Gibiansky and Sigmund \cite{gibsig} and  Lui \cite{lui} found new optimal multimaterial structures. In the past few years (2009-2012),  we suggested \cite{cherk09,cherkyuan,wheel,cherkgz} a new approach for optimal bounds of multimaterial mixtures and tested it on several examples of conducting composites.  

Here we consider  a  problem about two-dimensional multiphase  composites of a minimal compliance i.e. of a maximal stiffness.  Assume that a unit periodic square  cell $\Omega\subset R_2$ ($\|\Omega\|=1$) is subdivided into $N$ parts $\Omega_1, \ldots \Omega_N$ of given areas $m_i=\|\Omega_i\|$, $m_i>0$; these parts are filled with different elastic materials. 
For simplicity in notations, we consider here materials with zero Poisson ratio corresponding to a quadratic stress energy 
\begin{equation}
W_i(\kappa_i, \sigma)=\frac{1}{2} \kappa_i \mbox{Tr}(\sigma^2),
\end{equation} 
where stress stress  $\sigma$ describes an equilibrium, that is it is symmetric and divergencefree, 
\begin{equation}
\sigma=\sigma^T, ~\nabla \cdot \sigma=0
 \label{div-free}
 \end{equation} 
and $\kappa_i$ is a compliance of $i$th material. Below, we label the $i$th material itself as  $\kappa_i$.   The compliance of a composite is a piece-wise constant function
\begin{equation}
\kappa(x)=\sum_{i=1}^N \chi_i(x) \kappa_i.
\end{equation} 
 Assume that a given  external homogeneous stress $\sigma_0= \langle \sigma \rangle $ is applied to the periodic medium, where $\langle  \cdot \rangle$ means the average over a unit periodicity cell $\Omega$.
We look for the composite geometry that best adapts itself to the load $\sigma_0$ and we consider the following three closely related optimization problems: 
\newline
(i) Find the layout $\chi$ of materials in a periodicity cell, that  minimizes the energy $J(\chi, \sigma_0)$ of a composite:
\begin{equation}
I(m,\sigma_0)=\inf_{\chi_i \in {\cal M}} J(\chi, \sigma_0), \quad J(\chi, \sigma_0)=\inf_{\sigma\in {\cal U}}  \int_{\Omega_i }\chi_i W_i(\sigma) \,dx 
\label{1.1}
\end{equation}
where $\chi_i$ is the index function of $i$th material's domain,
$\chi_i(x)=1 $ if $x\in \Omega_i$ and $\chi_i(x)=0 $ if $x\not\in \Omega_i$.
\begin{eqnarray}
 {\cal U} &=& \left\{\sigma=\sigma^T, ~\nabla \cdot \sigma=0, \sigma_{ij}\in L_2(\Omega), ~\sigma \mbox{ is } \Omega-\mbox{periodic}, 
  \langle  \sigma \rangle = \sigma_0\right\}, 
  \label{1.2}
 \\
 \quad {\cal M} &=&  \left\{\chi_i: ~  \langle \chi_i \rangle=m_i\right\}, \quad m=(m_1, \ldots, m_n), \label{1.3}
 \end{eqnarray}
The problem (\ref{1.1}) of the stiffest three-materal elastic composite  in which the stress energy $W(\sigma)$ is minimized was studied in \cite{cherkgz}, the results are illustrated in Figure \ref{opt-str}. Similar results for  a electrical or thermal conducting composites were obtained earlier in \cite{cherkyuan} by a similar  approach. 
\\
(ii) Find the {\em G-closure}  \cite{lc82} -- the set of effective tensors $K_*(\chi)$ that correspond to all structures with fixed volume fractions $m$ of materials. A  boundary component of G-closure corresponds to effective compliance of extremal composite with the energy described in (\ref{1.1}), 
\begin{equation}
I(m, \sigma_0)=  \frac{1}{2} W(K_*(m,p),\sigma_0)
\end{equation} 
The optimal effective compliance tensor $K_*$ does not depend on the magnitude of $\sigma_0$, but depends on the 
 the ratio $p$ of eigenvalues of $\sigma_0$,  $K_*(m,p)$.  Varying this ratio and and computing the corresponding $K_*(m, p)$ we obtain the set of  optimal effective properties, that belongs to a component of G-closure boundary. 
 
{ \bf Remark} {\it In order to find the whole boundary of G-closure, one should minimize the sum of energy caused by several linearly independent excitations and vary parameters $p$ of these excitations and their magnitudes, see \cite{mybook,miltonbook} Such description was obtained in \cite{cherkyuan}  for  three-phase conducting composites.}
 %
\newline
(iii) Find  the quasiconvex envelope of multiwell Lagrangian. It is obtained from (\ref{1.1}) - (\ref{1.3}) if we  introduce  {\em costs} or  {\em weights} $\gamma_i$ for the unit of each material  instead of fixing their volume fractions and minimize the stress energy $J(\chi, \sigma_0)$  plus the cost of the used materials $\sum_i \gamma_i m_i$. The Lagrangian of this problem is
$$ F(\sigma, \gamma)= \inf_\chi \sum_{i=1}^N \chi_i \left( W_i(\sigma) + \gamma_i \right)
$$
We minimize over $\chi$ and obtain 
\begin{equation}
 F(\sigma, \gamma)= \min_{i=1,..., N}\left \{ W_i(\sigma)+ \gamma_i\right\} 
 \label{1.11}
 \end{equation} 
 where $\chi=( \chi_1, \ldots, \chi_N)$ is a nonconvex function of $\sigma$ and  $\gamma_i$. More specific, $F(\sigma,\gamma)$ is the minimum of several convex functions $W_i(\sigma)+ \gamma_i$ called {\em wells}. Notice that we identify the material (the well) by the value of minimizer $\sigma$.  A minimizing sequence  $\{\sigma^{(k)}\}$ oscillates and takes values in several wells in subdomains $\Omega_i$.
 
 The materials in an optimal composite are naturally ordered: larger values of $|\sigma|$ correspond to smaller values of $\kappa_i$. 
For some $\sigma_0$,  an optimal composite may degenerate into a two-material composite or a pure material. The dependence of cost $\gamma_i$ on the volume fraction $m_i$  is monotonic but may be not continuous \cite{bcd}.

 The minimizers in this nonquasiconvex problem oscillate in an infinitely  fine scale.
 The problem does not have a classical solution but only minimizing sequences. The dependence of the microstructures is due to the fact that the discontinuous stresses in neighboring grains keep the normal $n$ projection continuous, $[\sigma\cdot n ]^+_-$ which means that $\sigma$ depends on the geometry. 

\section{Relaxation and Quasiconvexity}
\paragraph{Quasiconvex envelope} 
The variational problem with nonquasiconvex Lagrangian (\ref{1.11}) should be be relaxed, which means that the oscillating sequences are to be replaced by effective (average over the periodicity cell) stresses in optimal composites; the properties of these composites also vary  from point to point in response to varying applied  stress $\sigma_0$, but this variation is slow and it can be neglected in the microscale when the optimal structure is determined. 
 
 The relaxed Lagrangian is defined by the {\em quasiconvex envelope} $QF$ that represents the energy of an optimal microstructure (more exactly, a limit of a sequence of microstructures) plus the cost. 
 \begin{eqnarray} QF(\sigma, \gamma) = \inf_{\zeta\in Z} \int_\Omega F(\sigma+ \zeta, \gamma) \, dx
\label{qf}\\
Z =\left \{\zeta:~\int_O \zeta=0, ~\nabla\cdot \zeta=0, ~ \zeta=\zeta^T, ~\zeta \mbox{ is }\Omega- \mbox{perodic} \right\}  \label{dc}
\end{eqnarray}
  where  $\gamma=(\gamma_1, \ldots, \gamma_n)$ is the deviation of the stress field from its average value $\sigma$. Notice that without the differential constraint   $ \nabla\cdot \zeta=0 $, (\ref{qf}) defines a convex envelope of nonconvex function $F(\sigma)$.  If the effective optimal compliance tensor $K_*$ is known,  $QF(\sigma, \gamma)$  can be conveniently expressed through it:
\begin{equation} QF(\sigma, \gamma) = \min_m \left(\frac{1}{2} W(K_*(m),\sigma_0) +m^T \gamma\right) .
\label{1.4}\label{gamma}
\end{equation}

Quasiconvexity was intensively studied in the last two decades, see for example \cite{dacorogna,mybook} and references therein; however, there are only a few examples of explicitly constructed multiwell quasiconvex envelope (Four gradients \cite{4-grad}, special case of Hashin-Shtrikman bounds \cite{gibsig,milton82,ls-multi,miltonbook}) which show that the technique for such problems is not yet developed. The first example of a component of quasiconvex envelope for a three-well Lagrangian is demonstrated below in Section \ref{inlarge}.

{\bf Remark on multi-face convex envelope } {\it The variety of microstructures in Figure \ref{opt-str} reveals the geometric complexity of the quasiconvex envelope that is a multiface  surface in the space of eigenvalues of $\sigma$. To illustrate this complexity of multiwell quasiconvex envelope, consider a simpler construction of convex envelope of minimum of several paraboloid wells.  The number $N_c$ of supporting points of each point of the envelope ${\cal C} W$ is defined by the number $N$ of strictly convex wells (no more than one point in a well) and by the dimension of the minimizer $d$ (Caratheodory theorem): $N_c\leq \min (N, d+1)$. For a two-well problem, the number of supporting points of the envelope is always two:   
The convex envelope of the minimum of two paraboloids in $R^n$ is either the paraboloids themselves or a cone stretched on them.  
In contrast,  the convex envelope of minimum of three arbitrary located paraboloids in $R^n$ consists of a flat component supported by all three paraboloids,  parts of conical surfaces supported by pairs of the paraboloids, or the paraboloids themselves. This geometric complexity is addressed in bounds for multimaterial composites. }

\section{Bounds of multimaterial composite properties} \label{newbound}
 Physically, $QF(\sigma_0)$ is the energy of the best composite plus the cost of the materials used; it is function of the average field $\sigma_0$. The relaxation (calculation of $QF$) is performed by a two-step procedure: (i) finding the exact  lower bound for the quasiconvex envelope and (ii) approximating these bounds by computing the energy in a class of microstructures similar to those shown  in Figures \ref{opt-str} and \ref{wheels}.  Calculating the lower bounds, we determine sufficient conditions on the fields in materials (wells) in an optimal structure \cite{mybook,miltonbook,cherk09}.   These also provide a hint for optimal structures such as high-rank laminates \cite{mybook,miltonbook} or wheel-type structures \cite{wheel}.


\paragraph{Principles of bounds derivation}
Deriving the bounds, we keep in mind the basic features of optimal layouts.

(i) Differential constraint $\nabla\cdot \sigma=0 $ in (\ref{1.2})) cannot be directly applied for solutions in domains of uncertain geometry. A lower bound for the problem (\ref{1.1}) - (\ref{1.3}) requires that  these  constraints are weakened and replaced by either integral or pointwise constaints. 

(ii) The energy of a a composite depends on the properties $\kappa=[\kappa_1, \ldots, \kappa_N]$ of materials, their volume fractions $m=[m_1, \ldots, m_N]$, the applied load $\sigma_0$, and the geometry of the structure. The bound is an infimum of the energy over all possible geometries, therefore it is a function of the first parameters only. The stress field is also function of the same parameters; moreover, the stress tensor in an optimal structure is independent of the distance from the boundary and other geometrical parameters, because all geometries are compared. Essentially, there is no difference between interior and boundary points, because the boundary could be arbitrary close to each point. This implies that in optimal structures some invariants of stress tensor are constant in each material. 

\paragraph{Translation bound}

Without  the constrain $\nabla\cdot \sigma=0 $, the variational problem (\ref{1.2})) is reduced to the finite-dimensional minimization problem. The lower bound of energy is represented by  the convex envelope of Lagrangian $F$. The effective compliance is  is computed as   the so-called Wiener bound  This bound is rough or not achievable. 

The Translation method  \cite{lc82,tartar,lc87,mybook,miltonbook} replaces the differential constraints by an integral constraint    
\begin{equation}
\langle \det \sigma \rangle= \det \langle  \sigma \rangle
\label{translator}
\end{equation}
 that follows from the constraint $\nabla \cdot \sigma=0$ and the Green's formula, see for example \cite{mybook}. 
  In the Translation method, equality (\ref{translator}) (the quadratic form of stress components) is added with a Lagrange multiplier $t$ to the finite-dimensional minimization problem that defines the convex envelope; the problem becomes a convex envelope of the {\em translated} wells $W_T= W_i(\kappa_i, \sigma)- t \det \sigma$. In order to obtain a proper envelope (not equal to $-\infty$), the translated quadratic wells $W_T$ should remain nonnegative, which constrains the range of $t$. 
  
  The obtained bound  is proven to be {\em exact} in many examples for two-well energy \cite{mybook,miltonbook}. However, the bound becomes rough for multiphase problem. 
For example, consider the lower Hashin-Shtrikman bound $\kappa_{HS}$ that is a special case of the translation bound for effective  compliance $\kappa_*$ of an isotropic composite (with zero Poisson ratio)
\begin{equation}
\kappa_* \ge \kappa_{HS}= - \kappa_1+ \sum_i^3\frac{m_i}{\kappa_i+\kappa_1}, \quad \mbox{if } \kappa_1< \kappa_2, \kappa_3.
\label{HS}
\end{equation}
 An addition of a conducting material $\kappa_1$ with  zero [sic!] volume fraction $m_1=0$ to a two-component isotropic composite of materials $\kappa_2$ and $\kappa_3$ surprisingly changes it because the bound depends on $\kappa_1$ even in the limit $m_1= 0$. This shows that bound is not exact for small $m_1$.

\paragraph{Pointwise constraints and supporting points} The multiwell energy optimization requires an account for other constraints besides (\ref{translator}) that also follow from differential properties (\ref{div-free}) of minimizer and from optimality requirements. 
Stress tensors in the materials in an optimal composite are constrained.  When these constraints are added to the translation method scheme, they result in a tighter bound. For several known examples, the new constraints  produce an exact bound for a multiwell Lagrangian \cite{cherk09,cherkyuan,wheel,cherkgz}.  The constraints are also satisfied for two-well Lagrangians as well, but there they do not change the bound and  they were not noticed. 
\paragraph{Equilibrium requirement \cite{cherk09}} 
We call the values of stress tensors in an optimal composite {\em supporting points} in the well. Physically, the supporting points are the alternating values of optimal stresses in the components of the structure (or limits of these values).  The quasiconvex envelope $QF(\sigma_0)$ (energy of an optimal composite) is the limit of the averaged energy computed on these oscillating minimizers.  All supporting points $\rho$ are located on the boundary of regions where quasiconvex envelope $QF(\sigma)$ coincides with Lagrangian $F(\sigma)$; at the points $\sigma= \rho$ of support,  the graph of the quasiconvex envelope $QF$ touches the graph of Lagrangian $F(\sigma)$.   We denote by $R_i$ the set of supporting points in the $i$th well (this set may be empty). The minimizers $\sigma^{(k)}(x)$ oscillates between these values: $\sigma^{(k)}(x) \to \rho \in  R_i$  if $x\in \Omega_i$.

By virtue of equilibrium, the normal stress is continuous at boundaries of grains in microstructures, which implies that Rank\,$ [\sigma]^+_-=1$ at any boundary.  However, this condition does not imply that all supports are in rank-one connection, the continuity may apply to the  stress tensor that is averaged in a smaller scale in the exterior of a grain. The normal supporting stress $\rho_\alpha \cdot n $ ($n$ is a normal) in a neighborhood of a boundary point inside a material grain is equal to the normal stress in the exterior neighborhood of this point where  another material or a smaller-scale mixture of other materials is located. Because the external stress in each material is a supporting stress $\rho$, their mixture belongs to the convex envelope ${\cal C} R_{A} $ of their support sets $R_{A}$. Each supporting point $\rho_\alpha $ of the problem (\ref{1.1}) is therefore in a rank-one contact  with the stress tensor $\rho_\alpha $ belonging to the convex envelope $ {\cal C} R_\alpha$ of the other supporting points, $ R_\alpha=R -\rho_\alpha$:

\begin{equation}
\exists \, \rho_A\in {\cal C} R_\alpha: ~\mbox{Rank}\,(\rho_A-\rho_\alpha)=1 ~~\Rightarrow ~~
 \rho_\alpha\in \left[\min_{\sigma\in R_A} \lambda_1(\sigma), ~ \max_{\sigma\in R_A}  \lambda_2(\sigma) \right] 
\label{compatibility}
\end{equation}
where $ \lambda_1$ and $\lambda_2$, $\lambda_1 \leq \lambda_2$ are  the minimal and maximal eigenvalues of tensor $\sigma$, respectively. Two {\it corollaries}  follow:
\\ 
1. If one of the wells (void) is supported by a single point  $\rho_\alpha=0$ (stress in void is zero), then 
the convex envelope ${\cal C} R_A$  of other supports must contain at least one nonpositive defined stress $\rho_A\in {\cal C} R_A$, such that   $\det(\rho_A-\rho_\alpha)= \det(\rho_A)=0$. 
\\
2. Assume that  all wells but one are supported by single isotropic minimizers $\rho_i= \beta_i I $ (stress in all material but the first one is constant), $\beta_2> \cdots > \beta_n$, and the first well is supported by a part of the line $\lambda_1+ \lambda_2= 2 \beta_1 $.   These conditions describe the fields obtained by Translation method without constraints. Then only the  supports in the interval $ \beta_n< \lambda_1,  \lambda_2 <\beta_2$ satisfy (\ref{compatibility}) and only they are is compatible. 
Notice, that  (\ref{compatibility})  is only a necessary condition for $s_\alpha $ to be a supporting point, since not all points of ${\cal C} R_A$ may correspond to a structure (some points in ${\cal C} R_A$ may be incompatible).

\subparagraph{Mean field inequality \cite{cd14}} 
This condition compares the supports with the average stress  $\sigma_0$.  In an optimal composite, all supporting fields $\rho_1\in R_1$ in the material with the lowest $\kappa_1< \kappa_2 $  satisfy the inequality
\begin{equation}
\det(\rho_1-\sigma_0) \le 0, \quad \forall \rho_1\in R_1
\label{mfineq}
\end{equation}
The proof is based on a special structural variation (the interchange of two elliptical inclusions of optimal shapes, see \cite{mybook,cherkKucuk}). One can show that such variation decreases the cost of the problem if (\ref{mfineq}) is not satisfied. This inequality is proven for a problem of optimal 2d composites from several isotropic components, but it should be generalized to 3d case. 

\noindent
{\bf Remark} {\it The inequality explains, in particular, why the best material in an optimal structure tends to form an envelope around the core of other materials or substructures: The normal stress in the  outer layer is equal to the average stress. In Hashin-Shtrikman coated circles, the normal stress in outer annulus increases and tangential stress decreases; equality in (\ref{mfineq}) is achieved at the external radius.}

\paragraph{New lower bounds and optimal structures}
The outlined technique allows for finding optimal multimaterial bounds. Here we describe the results following the simplified version of \cite{cherkgz} (one-constant elasticity, eqs. (\ref{1.1})-(\ref{1.3})). A similar technique leads to exact bonds for multiphase conducting composites \cite{cherk09,cherkyuan,wheel,cherkgz}.

Consider  problem (\ref{1.1}): The inequality (\ref{compatibility})  states that $\det \sigma \geq 0 $  everywhere in $\Omega$ if $\det \sigma_0 \geq 0 $ (this also agrees with the Alessandrini-Nesi inequality \cite{Alessandrini-Nesi}).
Applying the translation method, inequality (\ref{compatibility}), and assuming that $\det \sigma_0 \geq 0 $, we redefine the
translated energy as follows:  if the inequality $\det \sigma \geq 0$ is violated, the energy is equal to $+\infty$ :
\begin{equation} V_i(\sigma, t)=\left\{ \begin{array}{ll}
                 \frac{1}{2} \kappa_i \mbox{Tr}\, \sigma^2- t \det \sigma & \mbox{ if } \det(\sigma) \geq 0 \\
+ \infty & \mbox{ if } \det(\sigma) < 0
                \end{array}
  \right. .
\label{local}
\end{equation}
Observe that $V_i(\sigma, t)$  grows quadratically  with $\|\sigma\| $ for
all values $t>0$ of translation parameter $t$, even if it becomes a nonconvex function of $\sigma $ when
$ t>\kappa_i$, which eliminates the paradox (\ref{HS}) of Hashin- Shtrikman bound that was caused by constraint $t\le \kappa_1$.  The modified  bound $P_L V_0(\sigma_0)$ of the multiwell problem is:
\begin{eqnarray}
P_L V_0(\sigma_0) &=& \max_{t\geq 0} L(\sigma_0, t), ~~L(\sigma_0, t)= \min_{\sigma \in {\cal E}:} \sum_n \int_\Omega \chi_n V_n \,dx + t \det \sigma_0  \label{bound}
\\
{\cal E} &=& \left\{ \sigma:  \sigma=\sigma^T, ~\int_\Omega \sigma\, dx =\sigma_0, ~(\ref{mfineq}) \mbox{ is satisfied} \right\}
\end{eqnarray}

The
structure of minimizing sequences  depends on whether or not the wells are convex. Minimizers always oscillate between the wells. 
In a nonconvex well $V_i$, minimizers also oscillate between  boundaries  defined by inequalities (\ref{compatibility}) and (\ref{mfineq}).  
With these adjustments, we first find  the sets $S_i$ of supports  ($\sigma \in S_i$ if $x \in \Omega_i$) and then apply  to (\ref{bound}) the translation bound technique. Eliminating the differential constraints, we define the lower  bound as a solution to the finite-dimensional problem of constrained optimization; the new bound follows. The voluminous expression for the anisotropic bounds can be found in \cite{cherk09,cherkyuan,cherkgz}.

\paragraph{Example: Bounds for  isotropic three-material composite \cite{cherk09}} In the  case  of $\kappa_1 < \kappa_2< \kappa_3=\infty$, 
 the isotropic effective compliance $\kappa_*$, it is bounded by simple inequalities
  \begin{equation}
 \kappa_*(m_1, m_2)  \geq \left\{\begin{array}{ll}
 -\kappa_1 + \left(\frac{m_1}{2 \kappa_1}+\frac{m_2}{\kappa_1+\kappa_2}
\right)^{-1} 
   &\mbox{ if }
 m_{11}\leq m_1 \leq 1,
\\
\kappa_2+2\frac{\kappa_1}{m_1}(1-\sqrt{m_2})^2   &\mbox{ if } m_{11}\leq m_1 \leq m_{12},
\\
-\kappa_2 + \left(\frac{m_1}{2 \kappa_1}+\frac{m_2}{2\kappa_2}
\right)^{-1}.
    & \mbox{ if } ~~0 \leq m_1\leq m_{12}.
 \end{array} \right.
\ \label{b2}
 \end{equation}
(notice the irrational dependence of $m_2$),  where the threshold values $m_{11}$ and $m_{12}$ are
 \begin{equation}
m_{11}=\frac{2 \kappa_1}{\kappa_2+\kappa_1}(\sqrt{m_2}-m_2), \quad
m_{12}=\frac{\kappa_1}{\kappa_2}(\sqrt{m_2}-m_2).
 \label{b6}
  \end{equation}
 Bound (\ref{b2}) and its anisotropic generalizations are exact; they is realized either by the assemblages shown in Figure 2 \cite{wheel} or by isotropic structures in the upper line of Figure 1 \cite{cherk09}.   
 The first line in (\ref{b2}) is Hashin-Shtrikman bound (\ref{HS}): this bound is exact for sufficiently large $m_1$.

 \paragraph{Structure of the quasiconvex envelope in anisotropic case}
The above inequalities and the translation method allows us to find the bounds in optimal structures also for anisotropic composites shown in Figure \ref{opt-str}. We outline the results \cite{cherkgz,cd14}, see also \cite{cherkyuan}.
 The bounds depend on volume fractions and compliance of materials and anisotropy parameter $p$ of $\sigma_0$, which define the translation parameter $t$.  The bounds assume various forms depending on  the above inequalities. They  are satisfied as equalities (become active) in different regions, see Figure \ref{opt-str}, right field.   Namely:  
   
   \noindent  --   In region A,  $t=\kappa_2$,  inequality (\ref{compatibility}) is active in the first well $V_1$.
   
    \noindent   --  In region B, $t\in (\kappa_1, \kappa_2)$,  inequality (\ref{compatibility}) is active in the first well $V_1$.
   
     \noindent --   In region C, $t\in (0, \kappa_2)$, inequality (\ref{compatibility}) is active and all fields are constant. 
   
     \noindent --   In region D,  $t=\kappa_1$, the bound becomes the classical translation bound.
   
  \noindent --   In region E,  $t\in (\kappa_1, \kappa_2)$, both inequalities (\ref{compatibility})  and (\ref{mfineq}) are active in complimentary parts of $V_1$.  
   
Optimal structures shown in Figures 1 and 2 are obtained using  the described sufficient conditions. We find the sets $R_i$ of supports  in all wells as a part of the process go deriving the bounds. To find an optimal structure,  we must determine a laminate geometry corresponding to these stresses, by enforcing the connectedness
(the neighboring layers may in turn be laminates themselves, then the compatibility is applied to the average field in the laminate), see \cite{can,cherkyuan,cherkgz}. The isotropic optimal wheel assemblages in Figure \ref{wheels} are obtained  \cite{wheel} using the effective field theory \cite{HS1}. A radius-dependent anisotropic laminate in the "spike" region is homogenized and effective properties of multicoated cycles are obtained by the separation of variables.

\section{Problem in large. Quasiconvex envelope} \label{inlarge}
To obtain the quasiconvex envelope, we solve problem (\ref{1.4}) calculating the energy $W(m_1,m_2, \sigma)= \frac{1}{2} \sigma K_*(m_1, m_2) \sigma $, where $K_*(m_1, m_2) $  is an effective compliance of the optimal composite, adding the cost of materials $\gamma_1$,  $\gamma_2$ and $\gamma_3$, respectively, and minimizing the sum with respect to volume fractions $m_1$ and  $m_2$. This calculation is performed for all types of optimal composites. If the bound is  explicitly known, so is the component of the quasiconvex envelope. Here we show the isotropic component  $ (\sigma=s I)$ of the quasiconvex envelope that can be  obtained by minimizing  energy (\ref{b2}) of an optimal isotropic structure plus its cost  with respect to $m_1$ and $m_2$.
\paragraph{Range of parameters}
Range of $\gamma$:  For simplicity in the notations, we normalize the costs, assuming that $\gamma_1=1$, $\gamma_2=\gamma$ and $\gamma_3=0$. The energy is an even function of $s$, and it is enough to consider the case  $s>0$. 

The three-material composites are optimal if 
\begin{equation}
\gamma \in (\gamma_a, \gamma_b), \quad \gamma_a=\frac{\kappa_2}{\kappa_1}, \quad \gamma_b=\frac{2\kappa_1}{\kappa_1+\kappa_2}.
\label{range}
\end{equation}
If $\gamma>\gamma_b$ ($\kappa_2$ is too expensive), only $\kappa_1$ and $\kappa_3$ are used in optimal compositions; optimal isotropic structures are Hashin-Strikman coated circles HS(13) ($\kappa_1$ is the envelope, $\kappa_3$ is the inclusion); when intensity $s$ is large, the pure strong $\kappa_1$ material is optimal. 

If $\gamma<\gamma_a$ ($\kappa_2$ is too cheap), the optimal structures are pure material or two-material composites. For small values of stress intensity $s$, coated circles HS(23) are optimal; when the intensity increases, pure material $\kappa_2$ becomes optimal, then  HS(12) circles become optimal, and then pure strong material $\kappa_1$ is optimal.  Notice that material $\kappa_2$ is an envelope in H(23) circles, but an inclusion in H(12) circles.

Consider now the case of intermediate $\gamma$. Depending on the interval of stress intensity $s$, we observe four regimes and four type of optimal structures:
\begin{equation}
\begin{array} {lll}U_1&:~ s\in\left[0,~  \rho_1\right] & W2 \mbox{ or } L(13,2,13) \nonumber \cr
U_2&:~ s\in\left[\rho_1,~\rho_2 \right] & \kappa_2 \nonumber \cr
U_3&:~ s\in\left[\rho_2, ~\rho_3\right] & HS(23)  \mbox{ or } L(12,1)  \nonumber \cr
U_4&:~ s\in\left[\rho_3,~\infty\right] & \kappa_1  \nonumber \cr
\end{array}
\label{gc1}
\end{equation}
where
$$
\rho_1=\frac{\gamma}{\sqrt{\kappa_1}} \quad \rho_2=   2\sqrt{ \frac{\kappa_1(1-\gamma)}{ \kappa_2^2- \kappa_1^2}}, \quad
\rho_3=\sqrt{ \frac{(1-\gamma)(\kappa_1+\kappa_2)}{\kappa_1 (\kappa_2-\kappa_1)}}
$$
Notice that the volume fractions of materials in optimal composites  in intervals $U_1$ and $U_3$ vary  depending on the stress intensity.

\paragraph{Optimal energy (quasiconvex envelope) in dependence of $s$} 
~

1. In the interval $U_1$ (small values of $s$)  optimal structures are L(13,2,13) or W2. The Lagrangian (see the middle line in (\ref{b2})) is
\begin{equation}
F_1(m_1, m_2, s)=\left[ \frac{ \kappa_2}{2}+\frac { \kappa_1\, \left( 1- \sqrt{m_2} \right) ^2}{ m_1} \right] s^2+m_2\gamma+ m_1
\label{qc1}
\end{equation}
Minimizing $F_1$ over $m_1$ and $m_2$, we find their optimal values 
\begin{equation}
m_1^{(1)}=\frac{\kappa_1}{\gamma} \left(\frac{\gamma}{\sqrt{\kappa_1}}- s\right)s,
\quad m_2^{(1)}=  \frac{\kappa_1}{\gamma^2}\,s^2
\label{qc2}
\end{equation}
and the optimal energy (quasiconvex envelope) $QF_1(s)=F_1\left(m_1^{(1)}, m_2^{(1)}, s\right)$
\begin{equation}
QF_1(s)=\frac{1}{2} \left(\kappa_2 - \frac{2 \kappa_1}{\gamma} \right) s^2+2\sqrt{\kappa_1} s 
 \quad \mbox{if } s \in  \left(0, ~\rho_1   \right)
\label{qc3}
\end{equation}
Notice that the coefficient by $s^2$ is negative, the Lagrangian is not convex. 
This regime is valid until $m_2^{(1)}$ (see \ref{qc2})) reaches the limit, $m_2^{(1)}=1$ at $s=s_{t1}$. At this point, the quasiconvex envelope touches  the second well.

2. In the interval $U_2$, pure intermediate material $\kappa_2$ is optimal:
\begin{equation}
QF_2(s)=\frac{1}{2} \kappa_2 s^2 + \gamma, \quad \mbox{if } s\in \left(\rho_1, ~\rho_2 \right)
\label{qc4}
\end{equation}

3. For larger $s$ in the interval $U_3$, the 1-2-second rank laminate L(12,1) or the equivalent Hashin-Shtrikman coated circles $HS(23)$ become optimal. Notice that in this regime $m_3=0$ and $m_2=1-m_1$.   The Lagrangian $QF_3(s)$ is
$$ QF_3(s)=\min_{m_1} F_3(m_1, s)
$$
where
\begin{equation}
F_3(m_1, s)=\left(-\frac{\kappa_1}{2} +\left(\frac{m_1}{2\kappa_1}+\frac{1-m_1}{\kappa_1+\kappa_2}\right)^{-1}\right)s^2+m_1+\gamma(1-m_1).
\label{qc5}
\end{equation}
 We find optimal value $m_1^{(3)}$ of $m_1$,
\begin{equation}
m_1^{(3)}= -2\frac{\kappa_1}{\kappa_2-\kappa_1}+s \sqrt{\frac{\kappa_1(\kappa_1+\kappa_2)}{(1-\gamma)(\kappa_2-\kappa_1)}}
\label{qc51}
\end{equation}
and $QF_3=F_3(m_1^{(3)}, s)$,
\begin{equation}
QF_3=-\frac{\kappa_1}{2}s^2 + 2\sqrt{\frac{ \kappa_1(\kappa_1+\kappa_2) (1-\gamma) }{\kappa_2-\kappa_1} } s +\frac{\gamma (\kappa_1+ \kappa_2)-2 \kappa_1}{\kappa_2-\kappa_1}.
\label{qc6}
\end{equation}
Here also the coefficient by $s^2$ is negative, the Lagrangian is not convex.

4. Finally, $\kappa_1$ is optimal for large stresses:
\begin{equation}
QF_4(s) =\frac{\kappa_1}{2}s^2 +1\quad s\in U_4.
\label{qc7}
\end{equation}

The isotropic strain $\epsilon =\frac{\partial W}{\partial \sigma} $ is a piece-wise affine nonmonotonic function of $\sigma$.  $\epsilon(\sigma)$ decreases in composite zones $U_1$ and $U_3$, and increases in the zones $U_2$ and $U_4$ of pure materials, as one can see from (\ref{qc3}), (\ref{qc4}), (\ref{qc6}), and (\ref{qc7}); this  reflects  nonconvexity of the quasinvex envelope. 
Notice that the strain depends on stress piecewise linearly and nonmonotonically, because the increase of stress caused the redistribution of materials in an optimal composite.  Notice also, that $\epsilon(0)=2\sqrt{\kappa_1}\neq 0$ because infinitesimal stress corresponds to an optimal composite with infinitesimal fractions of elastic materials; the product of compliance and the stress (the strain) is finite and not zero. 

\begin{figure}[htbp]
\begin{center}
\includegraphics[width=.7in,height=.7in]{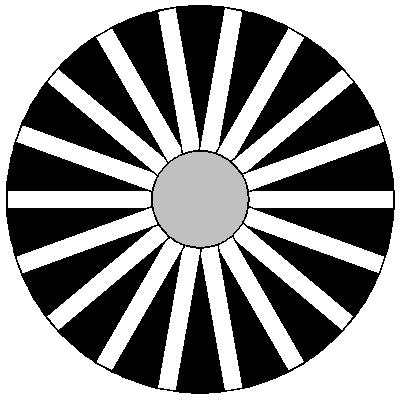}
\includegraphics[width=.7in,height=.7in]{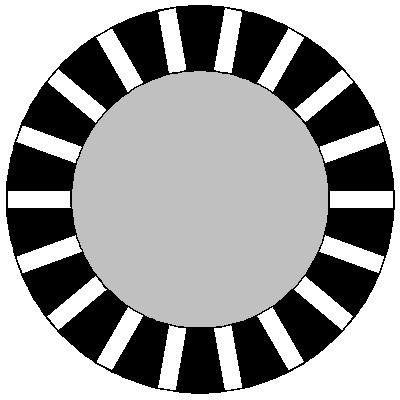}
~~~~~~~~\includegraphics[width=.7in,height=.7in]{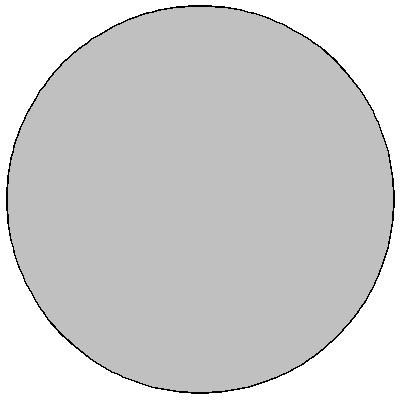} ~~~~~~
\includegraphics[width=.7in,height=.7in]{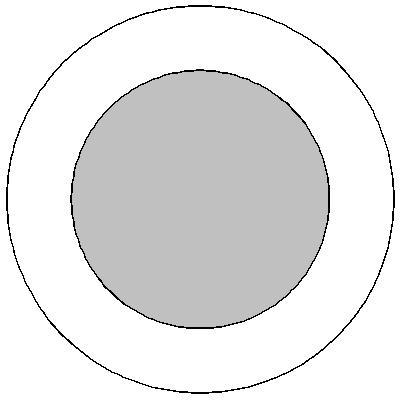} 
\includegraphics[width=.7in,height=.7in]{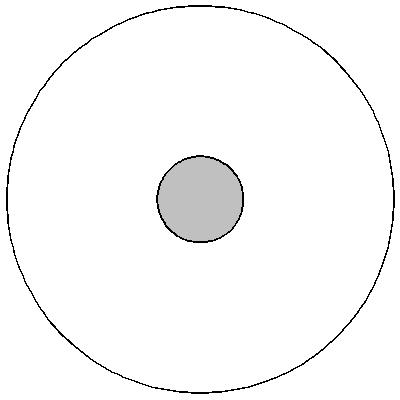}  ~~~~~~~~
\includegraphics[width=.7in,height=.7in]{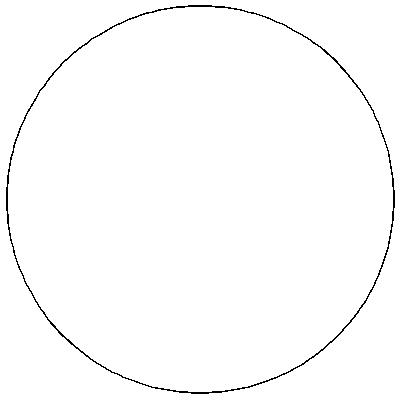} 
\caption{ Evolution of equivalent  isotropic optimal wheel assemblage with the increase of the magnitude of an applied hydrostatic stress: The variable three-material wheel, solid $\kappa_2$, the variable Hashin-Shtrikman coated spheres from $\kappa_2$ and $\kappa_1$, solid $\kappa_1$. 
}
\label{quasi}
\end{center}
\end{figure}

 In Figure \ref{quasi}, the evolution of {\em isotropic} wheel assemblages (Figure \ref{wheels}) is shown, the applied stress increases from left to right. When the stress is small, the structure consists of circular hubs of $\kappa_2$ jointed by strips of $\kappa_1$, between the strips is void $\kappa_3$. When stress increases, the  hubs grow and at a certain point the optimal composite becomes pure material $\kappa_2$. Next, it morphs to coated circles of $\kappa_2$ (inclusions) and $\kappa_1$ (envelope) and finally to a pure $\kappa_1$.

\paragraph{Anisotropic optimal structures}

\begin{figure}[htbp]
\begin{center}\includegraphics[width=2in,height=2in]{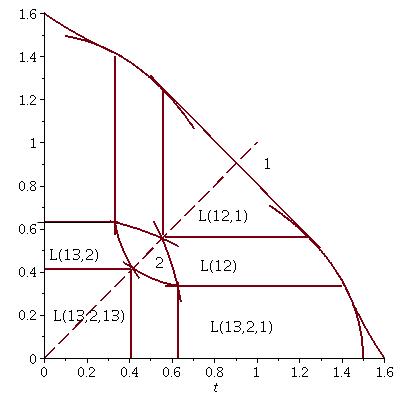} 
\caption{{\it Left:}  Types of optimal laminates (see Figure 1) in dependence on two eigenvalues of the external stress tensor, see \cite{bcd}. The case $\det \sigma >0$. The dashed line in the graph shows isotropic structures, see Fig.  \ref{quasi} 
}
\label{quasi1}
\end{center}
\end{figure}

 One can calculate the quasiconvex envelope of three-well energy for the general anisotropic stress as well. 
Without showing here the bulky formulas for the multifaced quasiconvex envelope, here we show  the evolution of optimal microstructures. Figure \ref{quasi1}   shows the placement of different types of  optimal microstructures, depending on the eigenvalues of the stress tensor. The shown case corresponds to the following range (\ref{range} of  cost $\gamma $.
Figure  \ref{quasi} shows zones 1 and 2 of pure first and second materials, the three-material composites shown in Figure \ref{opt-str}: L(13,2,13) are optimal for small $\|\sigma\|$, L(13,2)  and L(13,2,1) are optimal for anisotropic $\sigma$. Two-material laminates L(1,2) and second-rank laminates L(12,1) are optimal for larger $\|\sigma\|$. 
 Strongly anisotropic structures do not include zones of pure material $\kappa_2$. When only one eigenvalue of the applied stress increases, weak B-structures L(13,2,13) (see Figure 1) degenerate to C-structures L(13,2); these structures morph to E-structures L(13,2,1) and then to pure $\kappa_1$. 
 
Notice that not all structures in Figures \ref{opt-str} and \ref{wheels} are parts of the quasiconvex envelope for the chosen values of $\gamma$.  The remaining  structures  become components of the quasiconvex envelope   when $\gamma$ reaches the boundaries of the interval in (\ref{gamma}): A-structures correspond to $\gamma=\gamma_a$ and D-structures correspond to $\gamma =\gamma_b$. Outside this interval, no three-material structure is optimal.

\section{Structural optimization}

\begin{figure}[htbp]
\begin{center}
\includegraphics{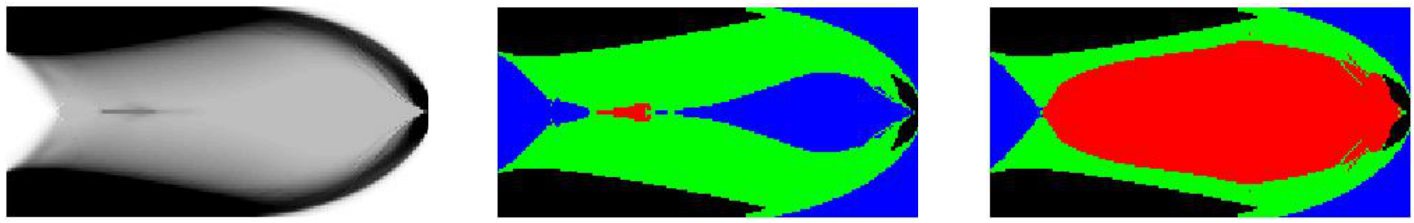}
\caption{Optimal two-materials  (gray) and multimaterial (colored) designs of a cantilever beam, from \cite{bcd}. The center field shows the case of low contrast of material properties, the right field shows  high contrast.  
In the colored designs,  the zone of the strongest material is shown in black, the zone of  weak composites L(13,2,13)  in blue, the zone of L(13,2,1) (strong, anisotropic composite) in green, and the zone of strong $(\kappa_1, \kappa_2)$-composites L(12,1) in red; compare with Figure \ref{quasi}. }
\label{opt-design}
\end{center}
\end{figure}

The most popular problem in structural optimization today, called ``topology optimization" \cite{sigmund}, is a problem of optimal layout of a material and void.  The optimal structural designs are commonly known as ``black-and-white" or ``grey" designs. Based on optimal multilateral composites, this suggested  approach allows us to instead deal with with ``multi-colored"  designs, see Figure \ref{opt-design}. 

 The first obtained optimal designs \cite{bcd} from three materials  are shown in Figure \ref{opt-design}. The designs are made from an expensive  strong material, a cheap weak material, and void.  The costs  are such that  $\gamma=\gamma_b$; in this case, zone 2 of pure $\kappa_2$ in Figure \ref{quasi} degenerates to a point. 

This design shows that the strong material tends to form elongated beam-like ligaments at the contour of the design while the concentration of the weaker material is larger in the inner areas of moderate and close to isotropy stress. For computations,  we adapted a  numerical algorithm of two-material structural optimization \cite{sigmund,grzegorz} based on steepest descent.

\paragraph{Conclusion}

At present, a theory of relaxation of multiwell Lagrangians  is not fully developed. The suggested methods for bounds and optimal microstructures work for special problems which are of independent interest for applications. A future development and generalization of these methods would allow for examination of a number of long-standing unsolved problems of optimal multilateral composites structures and will lead to generalization to multimaterial case of the obtained in the last 25 years results for optimal two-material composites.

\section*{Appendix: Exotic microstructures and metamaterials}

Structures with explicitly computable effective properties that are used in optimal bounds  play a special role in
the theory of composites. They allow for testing, optimizing, and demonstrating
the dependences of the structure and material properties, as well as for hierarchical
modeling of  complicated structures. They also permit for explicitly calculating 
fields inside the structure and tracking their dependence on structural parameters.
There are several known classes of such structures: laminates,  Hashin-Shtrikman coated
spheres structure \cite{HS1}, Schulgasser's structures \cite{schulgasser}, multiscale multi-coated spheres and multicoated laminates (see the discussion in \cite{lc-multi,gc-coupled}), or coated ellipsoids \cite{milton-benveniste}. 
These structures may or may not be optimal, but they all provide convenient and realistic models for the various sophisticated geometries that are used to create metamaterial hybrids between composites and lattices.  

\begin{figure}[htbp]
\begin{center}
\includegraphics[width=.9in,height=.9in]{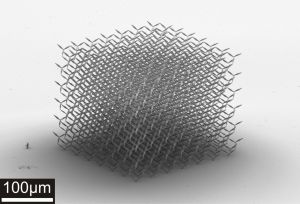} ~~~~~~
 \includegraphics[width=.9in,height=.9in]{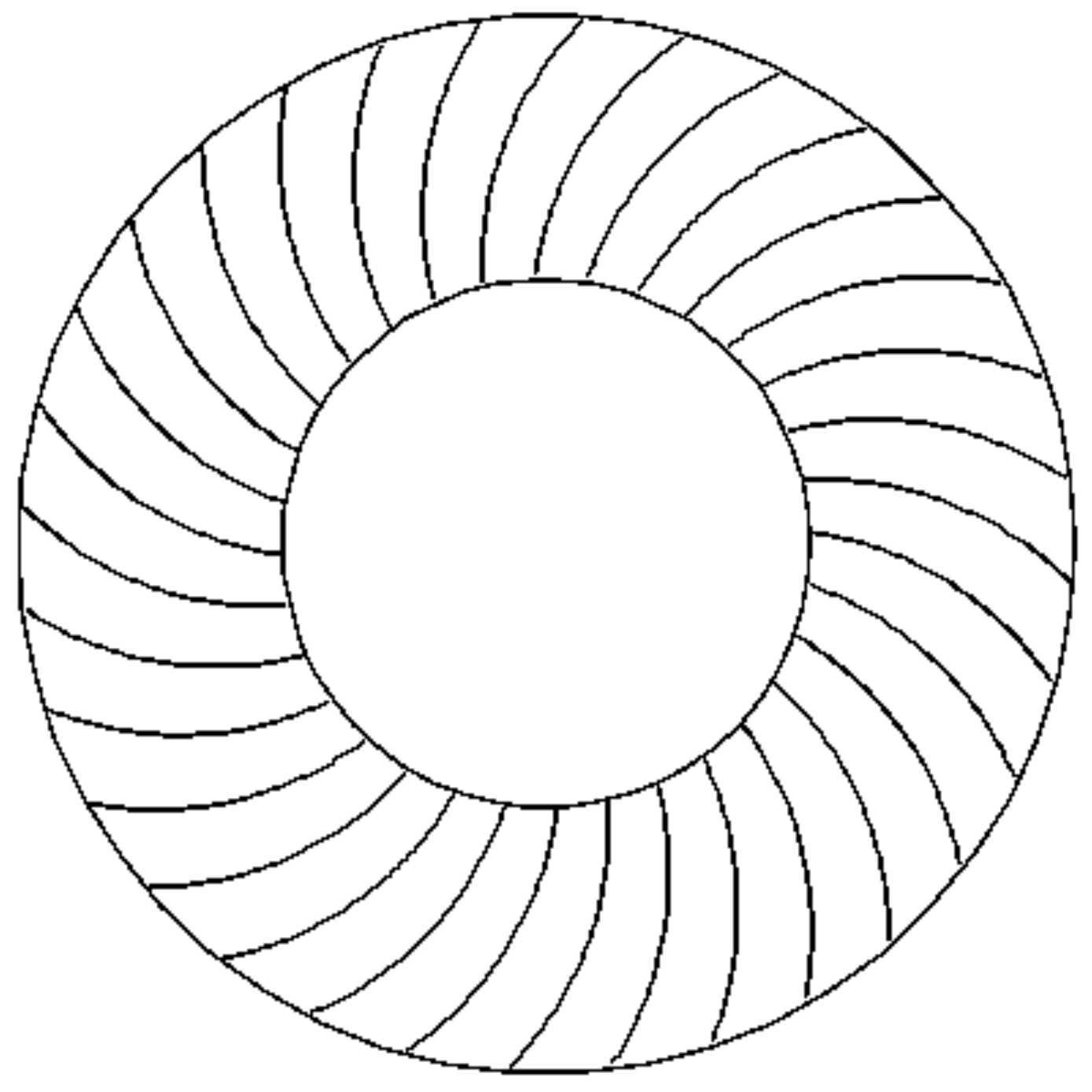} ~~ \includegraphics[width=.9in,height=.9in]{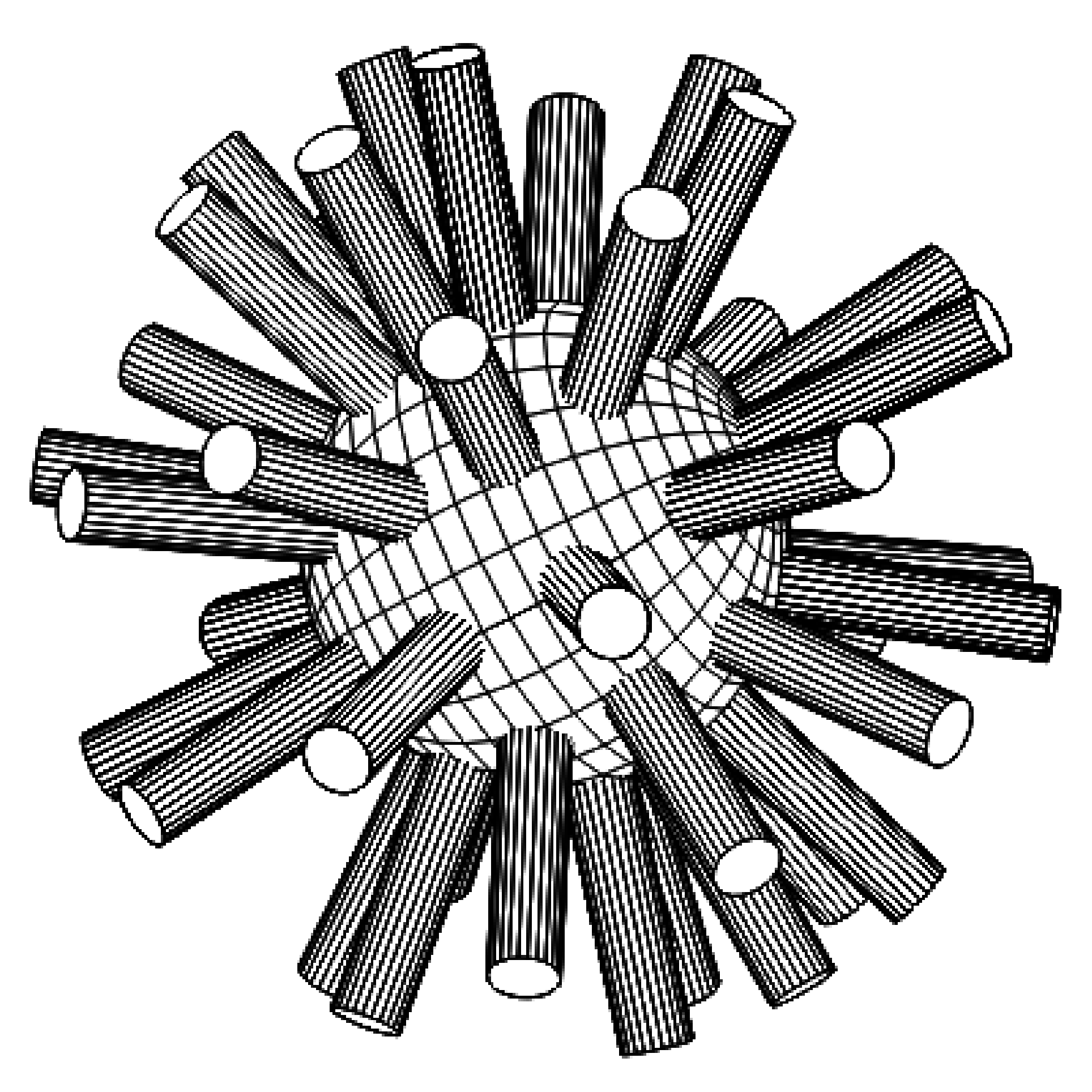} ~~ \includegraphics[width=.9in,height=.9in]{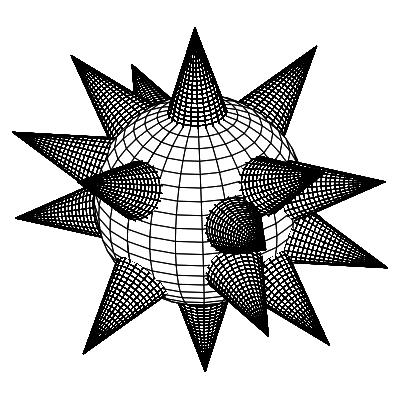} 
\caption{Left: Pentamode material that was experimentally produced by Marin Wegener et al. \cite{wegener} and was suggested in our theoretical paper \cite{miltoncherk95}. Center and right: Cartoons of exotic assemblage elements: Spiral with Core, Connected Hubs, and Spiky Balls, from  \cite{cherkpruss}. 
}
\label{exotic}
\end{center}
\end{figure}

One of the most exotic structures -- the pentamode -- that was suggested in our  paper \cite{miltoncherk95} in 1995 to prove the range of applicability of special classes of composites was constructed last year by the group of professor Martin Wegener  (KIT) \cite{wegener}, see Figure \ref{exotic}. Their experiment caught the attention of the mass media; 
 the material is promising for several industrial applications, such as an underwater acoustic invisibility cloak \cite{acoustics}.

Recently, we described new classes of such structures: the mentioned ``wheel assembly" \cite{wheel}, cylindrical assemblages of spirals with inner circular cylindrical
inclusions, spirals with shells and 3d assemblies which we call Connected Hubs and Spiky Balls, among others, see \cite{cherkpruss}.  The structures were investigated by the classical technique of Hashin-Shtrikman  coupled with hierarchical homogenization. These ``exotic" structures have interesting features of metamaterials.  For example, the spiral assemblies with inclusions (Figure \ref{exotic}) transform a homogeneous external current into a  homogeneous rotated current inside inclusions.  An observer there sees the current that flows, say, in a perpendicular direction to external  current ({\em "sun rises in the North"}).  The  spiral assemblages concentrate fields up to singularities in  central cores, which leads to an effective energy dissipation. These structures are natural metamaterials that may be used in electromagnetics and acoustics. In mechanics, these structures transform the overall pressure to a torque inside the cylindrical inclusions which should lead to interesting applications, such as sensors or compact electricity generators.  The Spiky Balls assemblages concentrate the current at the sharp edges, and Connected Hubs model a 3d network of connected reservoirs.
 
\paragraph{Acknowledgment} The author is thankful to Grzegorz Dzier\.{z}anowski and to Nathan Briggs for their comments and for providing graphics in Figure \ref{opt-design}.

\end{document}